\newcommand{\dataversione}{May~1st, 2025}

\documentclass[11pt,reqno,a4paper]{amsart}
\usepackage{amssymb}
\usepackage{mathrsfs}
\usepackage{fancyhdr}
\usepackage{hyperref}\hypersetup{colorlinks=true, citecolor=blue}
\usepackage{enumitem}
\usepackage[utf8]{inputenc}

\numberwithin{equation}{section}

\newtheoremstyle{mytheorem}
{10pt}
{10pt}
{\it}
{}
{\bf}
{.}
{ }
{\thmnumber{#2.~}\thmname{#1}\thmnote{~\rm#3}}

\newtheoremstyle{myremark}
{10pt}
{10pt}
{\rm}
{}
{\bf}
{}
{ }
{\thmnumber{#2.~}\thmname{#1}\thmnote{~\rm#3}}

\newtheoremstyle{myparagraph}
{10pt}
{10pt}
{\rm}
{}
{\bf}
{.}
{ }
{\thmnumber{#2.~}\thmname{#1}\thmnote{#3}}

\theoremstyle{mytheorem}
\newtheorem{theorem}[subsection]{Theorem}
\newtheorem{lemma}[subsection]{Lemma}
\newtheorem{corollary}[subsection]{Corollary}
\newtheorem{proposition}[subsection]{Proposition}

\theoremstyle{myremark}
\newtheorem{remark}[subsection]{Remark.}
\newtheorem{remarks}[subsection]{Remarks.\hskip-12pt}
\newtheorem{concludingremarks}[subsection]{Concluding remarks.\hskip-12pt}

\theoremstyle{myparagraph}
\newtheorem{parag}[subsection]{}
\newtheorem*{parag*}{}

\makeatletter
\def\@secnumfont{\sc}
\def\section{\@startsection%
{section}
{1}
\z@{1.5\linespacing\@plus .2\linespacing}
  {.7\linespacing}
  {\normalfont\sc\centering}}
\makeatother

\makeatletter
\renewenvironment{proof}[1][\proofname]{\par 
  \pushQED{\qed}%
  \normalfont \topsep10\p@\@plus6\p@\relax 
  \trivlist 
  \item[\hskip\labelsep 
    \bfseries 
    #1\@addpunct{.}]\ignorespaces 
}{%
  \popQED\endtrivlist\@endpefalse 
} 
\providecommand{\proofname}{Proof}
\makeatother

\newcommand{\footnoteb}[1]{\footnote{~#1}}
\newlist{enumeraterem}{enumerate}{1}
\setlist[enumeraterem]{label=(\roman*), leftmargin=0pt, itemsep=3pt, itemindent=30pt}

\newlist{enumeratethm}{enumerate}{1}
\setlist[enumeratethm]{label={\rm(\roman*)}, leftmargin=25pt, itemsep=2pt}

\newlist{itemizeb}{itemize}{1}
\setlist[itemizeb]{label=\textbullet, leftmargin=25pt, itemsep=2pt}

\newcommand{\R}{\mathbb{R}}
\newcommand{\N}{\mathbb{N}}
\newcommand{\Mass}{\mathbb{M}}
\newcommand{\Leb}{\mathscr{L}}
\newcommand{\D}{\mathscr{D}}
\newcommand{\Lip}{\mathrm{Lip}}
\newcommand{\Span}{\mathrm{span}}
\newcommand{\supp}{\mathrm{supp}}
\newcommand{\ii}{\mathbf{i}}
\newcommand{\eps}{\varepsilon}
\newcommand{\wrt}{with respect to\ }
\newcommand{\scalar}[2]{\langle #1 \, ; \, #2\rangle}
\newcommand{\bigscalar}[2]{\big\langle #1 \, ; \, #2\big\rangle}
\newcommand{\largewedge}{\mbox{\large$\wedge$}}
\DeclareMathOperator{\trace}{\mbox{\Large$\llcorner$}}

\begin{document}

    %
    %
\thispagestyle{empty}
~\vskip -1.1 cm

	%
	%
{\footnotesize\noindent 
[version: final, \dataversione]
\hfill J.~Funct.~Anal., 289 (2025), 111029 \par
\hfill DOI:~\href{https://doi.org/10.1016/j.jfa.2025.111029}{10.1016/j.jfa.2025.11102} \par
}

\vspace{1.7 cm}

	%
	%
{\centering\Large\bf
On the closability of differential operators
\\}

\vspace{.7 cm}

	%
	%
{\centering\sc 
Giovanni Alberti, David Bate, Andrea Marchese
\\}

\vspace{.8 cm}

	%
	%
{\rightskip .6 cm
\leftskip .6 cm
\parindent 0 pt
\footnotesize
{\sc Abstract.}
We discuss the closability of directional derivative operators with respect to a general Radon measure $\mu$ on $\smash{ \R^d }$; our main theorem completely characterizes the vectorfields for which the corresponding operator is closable from the space of Lipschitz functions $\Lip(\R^d)$ to $L^p(\mu)$, for $1\leq p\leq\infty$. 
We also discuss the closability of the same operators from $L^q(\mu)$ to $L^p(\mu)$, and give necessary and sufficient conditions for closability, but we do not have an exact characterization.

As a corollary we obtain that classical differential operators such as gradient, divergence and Jacobian determinant are closable from $L^q(\mu)$ to $L^p(\mu)$ only if $\mu$ is absolutely continuous with respect to the Lebesgue measure.

We finally consider the closability of  a certain class of multilinear differential operators;
these results are then rephrased in terms of metric currents.
\par
\medskip\noindent
{\sc Keywords:} 
closable operators, directional derivative operators, Lipschitz functions, Sobolev spaces, normal currents, metric currents.
\par
\medskip\noindent
{\sc MSC (2010):} 
26B05, 49Q15, 26A27.
\par

}

\section{Introduction}
\label{s1}
	%
	%
One way of defining the Sobolev spaces $W^{1,p}_0(\Omega)$ for an open set $\Omega$ in $\R^d$ is taking the completion of the space $C^1_c(\Omega)$ of functions of class $C^1$ with compact support on $\Omega$ with respect to the Sobolev norm $\|\cdot\|_{W^{1,p}}$.

This construction can be made more precise as follows:
we consider the graph of the gradient operator $\nabla: C^1_c(\Omega) \to C^0_c(\Omega;\R^d)$ as a subset of the product space $L^p(\Omega) \times L^p(\Omega;\R^d)$, take its closure $\Gamma$, and show that $\Gamma$ is still a graph, that is, for every $u\in L^p(\Omega)$ there exists \emph{at most one} $v\in L^p(\Omega;\R^d)$ such that $(u,v)\in\Gamma$.
We then consider the operator whose graph is $\Gamma$:
the domain is the Sobolev space $\smash{ W^{1,p}_0(\Omega) }$ and the operator is the gradient for Sobolev functions.%
\footnoteb{Moreover the norm of every $(u,v)\in\Gamma$ as element of $L^p(\Omega) \times L^p(\Omega;\R^d)$ agrees with $\|u\|_{W^{1,p}}$.}

\medskip
Note that the essential ingredient in this construction is that the closure of the graph of the gradient is still a graph.
The extension of this construction to more general operators leads to the following abstract definition:

\begin{parag*}[Closable operators]
Given $X, Y$ topological spaces, $D$ subset of $X$, and a map $T:D\to Y$, we denote by $\Gamma $ the closure of the graph
$\{(x,T(x))\colon x\in D\}$
in $X\times Y$, and we say that $T$ is \emph{closable} (from $X$ to $Y$) if $\Gamma$ is also a graph, that is, for every $x\in X$ there exists \emph{at most one} $y\in Y$ such that $(x,y)\in\Gamma$.
\end{parag*}

In this paper, we study the closability of certain first-order differential operators, and we focus in particular on directional derivative operators 
of the form \eqref{e:directional}.
The spaces $X$ and $Y$ are always spaces of functions on $\R^d$, taken from those listed in Paragraph~\ref{ss:functionspaces}.

Last but not least, we always intend continuity, closure and closability in the sequential sense. However, in many instances the sequential notions can be replaced by the topological ones for simple reasons (for example because the spaces are metrizable, or because we are dealing with a non-closability result).

\begin{parag}[Functions spaces]
\label{ss:functionspaces}
Through this paper $\mu$ is a Radon measure on $\R^d$ and
the space $Y$ is one of the following:
\begin{itemizeb}
\item
$L^p(\mu)$ with $1\le p\le\infty$, endowed with the strong topology;
\item
$L^p_w(\mu)$, which denotes the space $L^p(\mu)$ endowed with the weak topology if $p<\infty$, and with the weak* topology (as dual of $L^1(\mu)$) if $p=\infty$.
\end{itemizeb}
The space $X$ is one of the following:
\begin{itemizeb}
\item
$\Lip(\R^d)$, namely the space of Lipschitz functions on $\R^d$ endowed with the  following notion of convergence: $u_n\to u$ in $\Lip(\R^d)$ if $u_n\to u$ uniformly and the Lipschitz constants $\Lip(u_n)$ are uniformly bounded.%
\footnoteb{Equivalently, $u_n \to u$ uniformly and $\nabla u_n \to \nabla u$ in $L^\infty_w(\R^d)$. Therefore this notion of convergence is induced by a topology.}
\item
$L^q(\mu)$ or $L^q_w(\mu)$ with $1\le q\le\infty$.
\end{itemizeb}
Finally the set $D$ is always $C^1_c(\R^d)$.
\end{parag}

\begin{parag}[Directional derivative operators]
Let $v$ be a Borel vector field on $\R^d$;
we denote by $T_v$ the directional derivative operator on $C^1_c(\R^d)$ associated to $v$, that is, 
\begin{equation}
\label{e:directional}
T_v u  := \frac{\partial u}{\partial v} 
\quad
\text{for every $u\in C^1_c(\R^d)$.}
\end{equation}
\end{parag}

The next theorem is our main result.
The statement involves the notion of decomposability bundle $V(\mu,\cdot)$ of a measure $\mu$; the precise definition is given in Paragraph~\ref{ss:dec}, but for a first reading it suffices to know that $V(\mu,x)$ is a linear subspace of $\R^d$ for every $x\in\R^d$.

\begin{theorem}
\label{t:main1}
Let $v$ and $T_v$ be as above and assume that $v\in L^{p}(\mu)$ for some $1\le p \le \infty$.
\begin{enumeratethm}
\item
If $v(x)\in V(\mu,x)$ for $\mu$-a.e.~$x$, then every $u\in\Lip(\R^d)$ is differentiable at $\mu$-a.e.~$x\in \R^d$ in the direction $v(x)$, and the linear operator $\smash{ {\widetilde T}_v} : \Lip(\R^d)\to L^p_{w}(\mu)$ defined by
\begin{equation}
\label{e:extension}
\widetilde T_v u(x) := \frac{\partial u}{\partial v}(x)
\quad\text{for $\mu$-a.e.~$x\in\R^d$}
\end{equation}
is a continuous extension of $T_v$.
\\
It follows that $T_v$ is closable from $\Lip(\R^d)$ to $L^p_{w}(\mu)$. 

\item
Conversely, if $\mu(\{ x\colon v(x)\not\in V(\mu,x)\})>0$ then $T_v$ is nowhere continuous as an operator from $C^1_c(\R^d)$ (endowed with $\Lip$-convergence) to $L^p_{w}(\mu)$.
More precisely, for every $u\in C^1_c(\R^d)$ and every $\eps>0$ there exist a sequence $(u_n)$ in $C^1_c(\R^d)$ and $w\in L^p(\mu)$ with $w\ne T_v u$ such that
\begin{itemizeb}
\item
$u_n\to u$ uniformly;
\item
$\Lip(u_n)\leq\Lip(u)+\varepsilon$ for every $n$;
\item
$T_v u_n\to w$ in $L^p(\mu)$ if $p<\infty$, and in $L^\infty_w(\mu)$ if $p=\infty$.
\end{itemizeb}
\smallskip
It follows that $T_v$ is not closable from $\Lip(\R^d)$ to $L^p_w(\mu)$,
and not even  from $\Lip(\R^d)$ to $L^p(\mu)$ if $p<\infty$.
\end{enumeratethm}
\end{theorem}

\begin{remarks}
\label{r:aftermain}
\begin{enumeraterem}
%
\item 
Since the space $L^p_w(\mu)$ embeds continuously in $L^p(\mu)$,
the conclusion that $T_v$ is closable in Theorem~\ref{t:main1}(i) holds all the more so if we replace $L^p_w(\mu)$ by $L^p(\mu)$, but clearly $\widetilde{T}_v$ is not continuous from $\Lip(\R^d)$ in $L^p(\mu)$.

\item
The non-closability part of Theorem~\ref{t:main1}(ii) holds even if replace $\Lip(\R^d)$ by $L^q(\mu)$ with $1\leq q\leq\infty$ (or $L^q_w(\mu)$); this assertion is immediate when $\mu$ is a finite measure, because the fact that $u_n \to u$ uniformly implies that $u_n\to u$ in $L^q(\mu)$; if $\mu$ is only locally finite one should use that the sequence $(u_n)$ can be chosen so that the functions $u_n-u$ have uniformly bounded supports.
%
\end{enumeraterem}
\end{remarks}

We now turn our attention to the closability of classical differential operators such as gradient, divergence and Jacobian determinant.%
\footnoteb{The Jacobian determinant of $u\in C^1_c(\R^d;\R^d)$ is $Ju:=\det(\nabla u)$.}
Let indeed $T$ be any of these three operators:
it is well known that $T$ can be extended using a distributional definition to a continuous operator $\smash{{\widetilde T}}$ from $\Lip(\R^d)$ in $L^\infty_w(\Leb^d)$,%
\footnoteb{Here $\Leb^d$ is the Lebesgue measure, and, depending on which $T$ we consider, $\Lip(\R^d)$ and $L^\infty_w(\Leb^d)$ may denote spaces of $\R^d$-valued functions.}
and this implies that $T$ is closable from $\Lip(\R^d)$ to $L^\infty_w(\Leb^d)$.

It is natural to ask what happens if we replace the Lebesgue measure $\Leb^d$ by a general Radon measure $\mu$. The complete answer is contained in the following corollary of Theorem~\ref{t:main1}:

\begin{corollary}
\label{c:grad_div_jac}
Let $T$ be any of the following operators on $C^1_c(\R^d)$: gradient, divergence, Jacobian determinant, and let $1\le p\le\infty$. 
\begin{enumeratethm}
\item 
If $\mu$ is absolutely continuous \wrt the Lebesgue measure ($\mu\ll\Leb^d$)
then $T$ is closable from $\Lip(\R^d)$ to $L^p_w(\mu)$ and then also to $L^p(\mu)$.
\item 
If $\mu$ is not absolutely continuous \wrt the Lebesgue measure then 
$T$ is not closable from $\Lip(\R^d)$ to $L^p_w(\mu)$, and not even 
to $L^p(\mu)$ if $p<+\infty$.
\end{enumeratethm}
\end{corollary}

The next corollary answers a question posed by M.~Fukushima about the closability of the gradient (see \cite[Section~2.6]{Boga} and \cite{DLP}):

\begin{corollary}
\label{c:grad_Lp}
Let $T$ be as in Corollary~\ref{c:grad_div_jac} and take $1\le p<\infty$ and $1\le q\le\infty$.
Then $T$ is closable from $L^q(\mu)$ to $L^p(\mu)$ only if $\mu\ll\Leb^d$. 
\end{corollary}

\begin{parag*}[Structure of this paper]
In Section~\ref{s2} we collect some definition and preliminary results that are widely used through the rest of the paper, while Section~\ref{s3} is devoted to the proofs of the results stated above.

As pointed out in Remark~\ref{r:aftermain}(ii), Theorem~\ref{t:main1}(ii) gives a necessary condition for the closability of the directional derivative operator $T_v$ in \eqref{e:directional} from $L^p(\mu)$ to $L^p(\mu)$, but this condition is not sufficient (cf.~Remark~\ref{r:finalrem}(iii)).
In Section~\ref{s:Lp}, and more specifically in Theorem~\ref{t:lp_to_lp}, we give a sufficient condition for the closability of $T_v$; we do not know if this condition is also necessary.

In Section~\ref{ss:metric} we discuss the closability of a general class of alternating $k$-linear differential operators akin to the Jacobian determinant.
In Theorem~\ref{t:multilin} we give necessary and sufficient conditions for closability from $\Lip(\R^d)$ to $L^\infty_w(\mu)$; we do not know if these conditions match (unlike those in Theorem~\ref{t:main1}).
In the second part of Section~\ref{ss:metric} we reformulate these results in terms of metric currents in $\R^d$; among other things we obtain a reformulation of the \emph{Flat Chain Conjecture} (see \cite[Section~11]{metriccurrents}) both in terms of the $k$-tangent bundle of the measure associated to a metric current and in terms of closability/continuity of a suitably defined alternating $k$-linear differential operator (Theorem~\ref{r:flat_conjecture}).
\end{parag*}

\begin{concludingremarks}
\label{r:finalrem}
\begin{enumeraterem}
\item
Theorem~\ref{t:main1} has been used in \cite{BreGig} to give a new proof of the chain rule for $BV$ functions first proved in~\cite{AmbDal}; this new proof can be adapted to finite dimensional RCD spaces. 

\item
Through this paper, we consider for simplicity functional spaces defined on the domain $\R^d$. However, all results are essentially local in nature, thus the domain $\R^d$ can be easily replaced by any open subset of $\R^d$, and even by more general domains.

\item
Regarding Corollary~\ref{c:grad_Lp}, it is well known that the assumption $\mu\ll\Leb^d$ alone does not imply the closability of the gradient operator from $L^q(\mu)$ to $L^p(\mu)$, not even in dimension $d=1$.
For example, let $\mu$ be the restriction of the Lebesgue measure $\Leb^1$ to a totally disconnected compact subset of $\R$; then it is easy to prove that the derivative is not closable from $L^q(\mu)$ to $L^p(\mu)$ for any $1\le p, q\le\infty$.

A precise characterization of the (absolutely continuous) measures $\mu$ on $\R$ such that the derivative is closable from $L^2(\mu)$ to $L^2(\mu)$ has been given in \cite[Theorem~2.2]{Albev} (see also \cite[Theorem~3.1.6]{Fuk}).

\item
Theorem~\ref{t:main1} can be easily extended to more general \emph{first order} differential operators than just directional derivatives, thus obtaining a statement that includes Corollaries~\ref{c:grad_div_jac} and \ref{c:grad_Lp} as particular cases (see Remark~\ref{r:generalop}).

However, it seems that \emph{second order differential operators} are not easily included in our analysis. For instance, one would guess that the Laplace operator (defined on $C^2_c(\R^d)$) should be closable from $L^p(\mu)$ to $L^q(\mu)$ only if $\mu$ is absolutely continuous \wrt $\Leb^d$, but such statement does not seem to follow from any of the results in this paper (at least, not easily).

\item 
At the core of the proof of Theorem~\ref{t:main1}(ii), 
and of the ``only if'' part of Corollaries~\ref{c:grad_div_jac} and \ref{c:grad_Lp}, are the so-called ``width functions'', which are taken from  \cite[Lemma~4.12]{AlbMar}.
This notion was introduced by David Preiss while studying the differentiability of Lipschitz functions (cf.~\cite{ACP2}, \cite{Lindenstrauss2012}, and references therein); he also used it to give a first (unpublished) proof of Corollary~\ref{c:grad_div_jac}(ii) in dimension $d=2$.

\item
Given a Radon measure $\mu$ on $\R^d$, in \cite[Section~3]{BouChaJim} the authors define for $\mu$-a.e.~$x\in\R^d$ a tangent space $T_\mu(x)$ in such a way that the corresponding tangential gradient operator is closable from $\Lip(\R^d)$ to $L^\infty_w(\mu)$.
Theorem~\ref{t:main1}(ii) shows that $T_\mu(x) \subset V(\mu,x)$ for $\mu$-a.e.~$x$, and with some additional work one can prove that equality holds.
This remark gives a way to compute $T_\mu$ effectively; moreover, thanks to \cite[Theorem~1.1]{AlbMar}, every function $u\in\Lip(\R^d)$, 
is differentiable in the directions in $T_\mu(x)$ for $\mu$-a.e.~$x$.

\item 
Given a Radon measure $\mu$ on $\R^d$ and $1\le p\le\infty$,
in \cite[Section~2]{BouButSep} the authors define for $\mu$-a.e.~$x\in\R^d$ a tangent space $T^p_\mu(x)$ in such a way that the corresponding tangential gradient operator 
is closable from $L^p(\mu)$ to $L^p(\mu)$.
Using Theorem~\ref{t:main1}(ii) one can prove that  $T^p_\mu(x) \subset V(\mu,x)$ for $\mu$-a.e.~$x$. However, this inclusion may be strict  (cf.~the example in remark~(iii) above).
\end{enumeraterem}
\end{concludingremarks}

{\small
\begin{parag*}[Acknowledgements]
This research was initiated during visits of D.B.\ and A.M.\ at 
the Mathematics Department in Pisa. 
The visits of A.M.\ were partly supported by INdAM-GNAMPA.
The research of D.B.\ is supported by the European Union's Horizon 2020 Research and Innovation Programme, 
grant agreement no.~948021. 
The research of G.A.\ and A.M.\ is supported by INdAM-GNAMPA and by the Italian Ministry of University and Research via the project PRIN~2022PJ9EFL ``Geometric Measure Theory: Structure of Singular Measures, Regularity Theory and Applications in the Calculus of Variations'', funded by the European Union -- Next Generation EU, Mission~4, Component~2 -- CUP: E53D23005860006.
\end{parag*}
}

\section{Notation and preliminary results}
\label{s2}
Through this paper, sets and functions are always Borel;
measures are Borel and positive (unless stated otherwise) and, with the notable exception of Hausdorff measures, they are also locally bounded (that is, Radon).
If $\mu$ is a real- or vector-valued measure, $|\mu|$ denotes the variation of $\mu$.

\smallskip
In the next paragraph we recall the basic notation about currents in 
the Euclidean setting. For more details about currents see for instance \cite{KrantzParks}, \cite{Federer1996GeometricTheory}.

\begin{parag}[Classical currents]
A  \emph{$k$-dimensional current} $T$ in $\R^d$ is a continuous linear functional on the space $\D^k(\R^d)$ of smooth $k$-forms on $\R^d$ with compact support. 
The \emph{boundary} of $T$ is the $(k-1)$-current $\partial T$ defined by $\scalar{\partial T}{\omega} := \scalar{T}{d\omega}$ for every $\omega \in \D^{k-1}(\R^d)$.

The \emph{mass} of $T$, denoted by $\Mass(T)$, is the supremum of $\scalar{T}{\omega}$ over all $\omega\in\D^k(\R^d)$ such that $|\omega|\le 1$ everywhere. 
A current $T$ is called \emph{normal} if both $T$ and $\partial T$ have finite mass.

By Riesz theorem a $k$-current $T$ with finite mass can be viewed as a finite measure with values in the space $k$-vectors on $\R^d$, and therefore it can be written in the form $T=\tau\mu$ where $\mu$ is a finite positive measure and $\tau$ is a $k$-vector field in $L^1(\mu)$.
Thus the action of $T$ on a $k$-form $\omega$ is given by
\[
\scalar{T}{\omega} 
= \int_{\R^d} \scalar{\tau(x)}{\omega(x)} \, d\mu(x)
\, .
%
%
\]
Given a Lipschitz path $\gamma:[a,b]\to\R^d$, we denote by $[\gamma]$ the associated current, namely the $1$-dimensional current defined by
\[
\bigscalar{[\gamma]}{\omega} 
= \int_a^b \bigscalar{\gamma'(t)}{\omega(\gamma(t))} \, dt
\, .
\]
\end{parag}

\begin{parag}[Decomposability bundle]
\label{ss:dec}
Given a measure $\mu$ on $\R^d$, its \emph{decomposability bundle} is a Borel map $V(\mu,\cdot)$ on $\R^d$ whose values are linear subspaces of $\R^d$ defined as follows: a vector $v\in\R^d$ belongs to $V(\mu,x)$ if and only if there exists a $1$-dimensional normal current $N$ in $\R^d$ with $\partial N=0$ such that
\begin{equation}
\label{2.1}
\lim_{r\to 0}\frac{|N-v\mu|(B(x,r))}{\mu(B(x,r))}=0
\, .
\end{equation}
The decomposability bundle $V(\mu,x)$ was introduced in \cite{AlbMar}, Paragraph~2.6;
in Paragraph~6.1 of the same paper, the authors introduced also an auxiliary bundle $N(\mu,x)$, and then they proved that these two bundles agree \cite[Theorem~6.4]{AlbMar}.
The definition above is actually that of the auxiliary bundle $N(\mu,x)$, which is simpler to state than the original definition of the decomposability bundle.
\end{parag}

In the ensuing sections we will use several results from \cite{AlbMar} concerning the decomposability bundle. We state here the most relevant ones. 

\begin{theorem}[(see {\cite[Theorem~1.1]{AlbMar}})]
\label{t:rademacher4measures}
Let $\mu$ be a measure on $\R^d$. 
Then the following statements hold:
\begin{enumeratethm}
\item
Every Lipschitz function $f$ on $\R^d$ is differentiable at $\mu$-a.e.~$x$ with respect to the linear subspace $V(\mu,x)$, that is, there exists a linear function from $V(\mu,x)$ to $\R$, denoted by $d_V f(x)$, such that
\[
f(x+h) = f(x) + \scalar{d_V f(x)}{h} + o(|h|)
\quad\text{for $h\in V(\mu,x)$.}
\]

\item
The previous statement is optimal, meaning that there exists a Lipschitz function $f$ on $\R^d$ such that for $\mu$-a.e.~$x$ and every $v\notin V(\mu,x)$ the derivative of $f$ at $x$ in the direction $v$ does not exist.
\end{enumeratethm}
\end{theorem}

\begin{theorem}
\label{c:bundlesingular}
Let $\mu$ be a measure on $\R^d$. Then $V(\mu,x)=\R^d$ for $\mu$-a.e.~$x\in\R^d$ if and only if $\mu\ll\Leb^d$.
\end{theorem} 

\begin{proof}
The ``if'' implication is contained in \cite[Proposition~2.9(iii)]{AlbMar}.
The hard part is the ``only if'' implication,  which is a direct consequence of the results in \cite{AlbMar} and \cite{DPR};
for instance, it can be obtained by combining Theorem~\ref{t:rademacher4measures} and \cite[Theorem~1.14]{DPR}.
\end{proof}

\section{Proof of Theorem~\ref{t:main1} and Corollaries~\ref{c:grad_div_jac}, \ref{c:grad_Lp}}
\label{s3}
\begin{proof}[Proof of Theorem~\ref{t:main1}(i), case $\boldsymbol{p=\infty}$]
The operator ${\widetilde T}_v: \Lip(\R^d) \to L^\infty_w(\mu)$ given by formula \eqref{e:extension} is well defined thanks
to the assumption $v(x)\in V(\mu,x)$ for $\mu$-a.e.~$x$ and  Theorem~\ref{t:rademacher4measures}(i), and it is obviously an extension of 
the operator $T_v$ defined in \eqref{e:directional}.

It remains to prove that ${\widetilde T}_v$ is continuous.
By \cite[Theorem~6.3]{AlbMar} there exists a normal $1$-current $N=\tilde v \tilde\mu$ on $\R^d$ such that
\begin{itemizeb}
\item 
$\partial N=0$;
\item
$\tilde v \in L^{\infty}(\tilde\mu)$ and $\tilde v(x) \in V(\tilde\mu,x)$ for $\mu$-a.e.~$x$;
\item 
$\tilde v$ and $\tilde\mu$ extend $v$ and $\mu$ in the following sense: 
$\tilde\mu=\mu+\sigma$ with $\sigma$ and $\mu$ mutually singular, and $\tilde v(x)=v(x)$ for $\mu$-a.e.~$x$. 
\end{itemizeb}
Let $\smash{ {\widetilde T}_{\widetilde v} }: \Lip(\R^d) \to L^\infty_w(\tilde\mu)$ be the operator defined by formula \eqref{e:extension} with $v$ replaced by $\tilde v$;
one easily checks that the continuity of $\smash{{\widetilde T}_v}$ follows from that of $\smash{{\widetilde T}_{\widetilde v}}$.

\smallskip
To prove the continuity of $\smash{{\widetilde T}_{\widetilde v}}$, we note that for every $u\in\Lip(\R^d)$ the boundary of the current $uN$ is given by
\begin{equation}
\label{e:bordifN}
\partial(uN)= -{\widetilde T}_{\widetilde v} \, u \, \tilde \mu
\end{equation}
(use Proposition~5.13 in \cite{AlbMar} and the fact that $\partial N=0$) and then
\begin{equation}
\label{e:bordifN2}
\int_{\R^d} \varphi \, {\widetilde T}_{\tilde v} u \, d\tilde\mu
= -\scalar{uN}{d\varphi}
\quad\text{for every $\varphi\in\D(\R^d)$,}
\end{equation}
where $\D(\R^d)$ is the space of smooth test functions with compact support on $\R^d$.

\smallskip
Consider now a sequence $(u_n)$ such that $u_n \to u$ in $\Lip(\R^d)$. 
Since $u_n\to u$ uniformly, the currents $u_nN$ converge to $uN$ \wrt the mass norm $\Mass$, and therefore also in the sense of currents; then formula \eqref{e:bordifN2} implies that
\begin{equation}
\label{e:form_bordi}
\lim_{n\to\infty}
\int_{\R^d} \varphi \, {\widetilde T}_{\tilde v} u_n  \, d\tilde\mu
=
\int_{\R^d} \varphi \, {\widetilde T}_{\tilde v} u \, d\tilde\mu
\quad\text{for every $\varphi\in\D(\R^d)$.}
\end{equation}
Since $\D(\R^d)$ is dense in $L^1(\tilde\mu)$ and the functions $\smash{ {\widetilde T}_{\tilde v} } u_n$ are uniformly bounded in $L^\infty(\tilde\mu)$, \eqref{e:form_bordi} holds also for every $\varphi\in L^1(\tilde\mu)$, which means that
\[
{\widetilde T}_{\tilde v} u_n \to {\widetilde T}_{\tilde v} u
\quad\text{in $L^\infty_w(\tilde\mu)$,}
\]
and the continuity of ${\widetilde T}_{\tilde v}$ is proved.
\end{proof}

\begin{proof}[Proof of Theorem~\ref{t:main1}(i), case $\boldsymbol{p<\infty}$]
As for the case $p=\infty$, the operator ${\widetilde T}_v:\Lip(\R^d)\to L^p_w(\mu)$ is well defined and extends $T_v$.
To prove that ${\widetilde T}_v$ is continuous, we consider the vector field $\widehat{v}$ on $\R^d$ defined by
\begin{equation}
\label{e:renorm}
\widehat{v} (x):=
\begin{cases}
    \frac{v(x)}{|v(x)|} & \text{if $v(x) \ne 0$,} \\
    0 & \text{if $v(x) = 0$.}
\end{cases}
\end{equation}
Then $\widehat{v}$ is bounded and belongs to $V(\mu,x)$ for $\mu$-a.e.~$x$ because so does $v$;
hence $\smash{ \widetilde{T}_{\widehat v} }: \Lip(\R^d)\to L^\infty_w(\mu)$ is continuous by Theorem~\ref{t:main1}(i), case $p=\infty$.

Moreover the identity $v=|v| \widehat{v}$ implies that $\widetilde{T}_v u = |v| \, \widetilde{T}_{\widehat v} \, u$ for every $u$, and therefore, for every $\varphi\in L^q(\mu)$ where $q$ is the H\"older conjugate of the exponent~$p$, there holds
\[
\int_{\R^d} \big( \widetilde{T}_v u \big) \, \varphi \, d\mu
= \int_{\R^d} \big( \widetilde{T}_{\widehat v} \, u \big) \, |v| \varphi \, d\mu
\, .
\]
The continuity of $\widetilde{T}_{\widehat v}$ and the fact that $|v|\varphi$ belongs to $L^1(\mu)$ imply the continuity of the right-hand side (\wrt $u$), and therefore also of the left-hand side, which in turn implies the continuity of $\widetilde{T}_v$.
\end{proof}

\begin{proof}[Proof of Theorem~\ref{t:main1}(ii), case $\boldsymbol{p=\infty}$]
Given $y,y'\in\R^d$ and $W$ subspace of $\R^d$, we write
$\theta(y,y')$ for the angle between $y$ and $y'$, 
and $\theta(y,W)$ for the angle between $y$ and $W$ (both angles are set to be $0$ if any of the vectors involved is null).

Let $\vartheta(x):=\theta(v(x),V(\mu(x))$ for every $x\in\R^d$:
by assumption $\vartheta$ is strictly positive on a set of positive $\mu$-measure, and therefore we can find $\bar x$ such that $\vartheta(\bar x)>0$ and both  $\vartheta$ and $v$ are approximately continuous at $\bar x$.
Then, having set $\beta:=\frac{1}{3}\vartheta(\bar x)$ and $w:=v(\bar x)$, for every $\alpha$ with $0<\alpha<\beta$ the set
\[
E_\alpha := \big\{ x\in\R^d \colon 
   \vartheta(x) \ge 2\beta 
   \, , \ 
   \theta(v(x),w) \le\alpha
\big\}
\]
has positive $\mu$-measure. Moreover, for every $x\in E_\alpha$ there holds
\[
\theta(w,V(\mu,x)) 
\ge \theta(v(x),V(\mu,x)) - \theta(v(x),w) 
\ge 2\beta-\alpha
>\beta
\, , 
\]
which means that the intersection of $V(\mu,x)$ and the cone
$C:=\{y\colon \theta(w,y)\le\beta\}$ contains only $0$.
Therefore, by \cite[Lemma~7.5]{AlbMar}, $E_\alpha$ contains a compact set $F_\alpha$ with positive $\mu$-measure which is $C$-null in the sense of \cite[Paragraph~4.11]{AlbMar}, and then, by \cite[Lemma~4.12]{AlbMar}, there exists a sequence of smooth functions $f_n:\R^d\to\R$ such that, for every $x\in\R^d$,
\begin{itemizeb}
\item[(a)]
$0\leq f_n(x)\leq \frac{1}{n}$; 
\item[(b)]
$0\leq d_wf_n(x)\leq 1$, and $d_wf_n(x)= 1$ if $x\in F_\alpha$;
\item[(c)]
$|d_Wf_n(x)|\leq 1/\tan\beta$ where $d_W$ is the restriction of the differential $d$ to the subspace $W:=w^\perp$ (and $|\cdot|$ is the operator norm).
\end{itemizeb}
Using statements (b) and (c) we obtain that:
\begin{itemizeb}
\item[(d)]
there exists a finite constant $L$ such that $\Lip(f_n)\le L$ for every $n$;
\item[(e)]
for $\alpha$ small enough there exists a constant $\delta>0$  such that $T_v f_n(x) \ge\delta$ for every $x\in F_\alpha$ and every $n$.
\end{itemizeb}
We finally consider the functions $u_n := u + L^{-1}\varepsilon f_n$.
Statement~(a) implies that $u_n \to u$ uniformly; (d)~implies $\Lip(u_n) \le \Lip(u)+\eps$; (e)~implies $T_v u_n \ge T_v u + \delta'$ on $F_\alpha$, where $\delta' := L^{-1}\varepsilon\delta$, 
and therefore, possibly passing to a subsequence, we have that
$T_v u_n \to w$ in $L^\infty_w(\mu)$ for some $w\ne T_v u$.
\end{proof}

To prove Theorem~\ref{t:main1}(ii) for $p<\infty$ we need the following lemma.

\begin{lemma}
\label{p:mazur}
Let $1\leq p<\infty$ and $T:C^1_c(\R^d)\to L^p(\mu)$ be a linear operator.
Let $(u_n)$ be a sequence of functions in $C^1_c(\R^d)$ such that $u_n\to u$ uniformly and $Tu_n\to w$ in~$L^p_w(\mu)$, that is, weakly.
Then there exists a sequence $(\tilde u_n)$ of convex combinations of the elements of $(u_n)$ such that {\rm(a)}~$\tilde u_n\to u$ uniformly, and {\rm(b)}~$T\tilde u_n\to w$ in~$L^p(\mu)$, that is, strongly.    
\end{lemma}

\begin{proof}
Since $Tu_n \to w$ weakly in $L^p(\mu)$, 
by the version of Mazur's lemma stated in \cite[Lemma~10.19]{RenRog},
for every $n=1,2,\dots$ there exist an integer $N(n)\ge n$ and for every $k$ with  $n\le k \le N(n)$ there exist real numbers $\alpha^n_k \ge 0$ with sum equal to $1$ such that the functions
\[
w_n := \sum_{k=n}^{N(n)} \alpha^n_k \, Tu_k
\]
converge to $w$ strongly in $L^p(\mu)$. To conclude we set
\[
\tilde{u}_n:= \sum_{k=n}^{N(n)} \alpha^n_k\,  u_k
\, .
\]
Indeed (b) follows from the identity $T\tilde{u}_n=w_n$ and the fact that $w_n\to w$ in $L^p(\mu)$. Moreover, since $\tilde{u}_n -u$ is a convex combination of the functions $u_k -u$ with $k\ge n$, by the convexity of the supremum norm $\| \cdot\|$ there holds
\[
\| \tilde{u}_n - u\| \le \delta_n:=\sup_{k\ge n} \|u_k -u\|
\, ,
\]
and $\delta_n\to 0$ because $u_n\to u$ uniformly. Thus (a) is proved.
\end{proof}

\begin{proof}[Proof of Theorem~\ref{t:main1}(ii), case~$\boldsymbol{p<\infty}$]
Let $\widehat{v}$ be the vector field defined in \eqref{e:renorm}.
Then $\widehat{v}$ satisfies the assumptions of Theorem~\ref{t:main1}(ii), case~$p=\infty$, and we take $(u_n)$ to be the sequence given in that statement.

Using the identity $T_v u = |v| T_{\widehat v} \, u$ and the fact that  $T_{\widehat v} \, u_n\to w$ in $L^\infty_w(\mu)$ with $w\ne T_{\widehat v} \, u$, we obtain 
\[
T_v u_n = |v| \, T_{\widehat v} \, u_n \to w':=|v|w
\ \text{in $L^p_w(\mu)$, and $w'\ne T_v u$.}
\]
Finally we use Lemma~\ref{p:mazur} to construct from $(u_n)$ a new sequence $(\tilde u_n)$ such that $T_v u_n \to w'$ strongly in $L^p(\mu)$.
It is easy to check that the sequence $(\tilde u_n)$ satisfies all requirements.
\end{proof}

\begin{proof}[Proof of Corollary~\ref{c:grad_div_jac}(i)]
It is well known that each of the operators $T$ can be extended using a suitable distributional definition to a continuous operator from the Sobolev space $W^{1,\infty}(\R^d)$ to $\smash{ L^\infty_w(\Leb^d) }$; this means that $T$ can be extended to a continuous operator $\smash{ \widetilde{T}}: \Lip(\R^d) \to L^\infty_w(\Leb^d)$.

Since $\mu\ll\Leb^d$, the space $L^\infty_w(\Leb^d)$ embeds continuously in $L^\infty_w(\mu)$ and therefore $\smash{\widetilde{T}}$ is also a continuous operator from $\Lip(\R^d)$ to $L^\infty_w(\mu)$, and this implies that $T$ is closable from $\Lip(\R^d)$ to $L^\infty_w(\mu)$.

Moreover, $L^\infty_w(\mu)$ embeds continuously also in $(L^p_w(\mu))_{\mathrm{loc}}$,%
\footnoteb{And even in $L^p_w(\mu)$ if $\mu$ is a finite measure.}
thus $\smash{ \widetilde{T} }$ is also a continuous operator from $\Lip(\R^d)$ to $(L^p_w(\mu))_{\mathrm{loc}}$, and this implies that $T$ is closable from $\Lip(\R^d)$ to $(L^p_w(\mu))_{\mathrm{loc}}$.

Finally both $L^p_w(\mu)$ and $L^p(\mu)$ embed continuously in $(L^p_w(\mu))_{\mathrm{loc}}$, and this implies that $T$ is also closable from $\Lip(\R^d)$ to $L^p_w(\mu)$ and $L^p(\mu)$.
\end{proof}

\begin{proof}[Proof of Corollary~\ref{c:grad_div_jac}(ii)]
By the assumption on $\mu$, Theorem~\ref{c:bundlesingular} implies that $V(\mu,x) \ne\R^d$ on a set of positive $\mu$-measure. Then there exists an element $e_k$ of the canonical base of $\R^d$ such that $e_k\notin V(\mu,x)$ on a set of positive $\mu$-measure. 

Thus the constant vector field $e_k$ satisfies the assumption of Theorem~\ref{t:main1}(ii), and therefore there exist a sequence of functions $u_n\in C^1_c(\R^d)$ such that $u_n\to 0$ in $\Lip(\R^d)$ and the $k$-th partial derivatives $\partial_ku_n$ converge to some $w\ne 0$ in $L^p(\mu)$ 
if $p<\infty$, and in $L^\infty_w(\mu)$ if $p=\infty$.

To conclude we notice that the features of the sequence $(u_n)$ imply the non-closability of the gradient from $\Lip(\R^d)$ to $L^p_w(\mu)$, and even to $L^p(\mu)$ if $p<\infty$;
the sequence of vector fields $(u_n\, e_k)$ implies the non-closability of the divergence, and finally the sequence of maps $(U_n)$ given by $U_n(x) := x +u_n(x)\, e_k$ implies the non-closability of the Jacobian determinant.
\end{proof}

\begin{proof}[Proof of Corollary~\ref{c:grad_Lp}]
This statement follows from Corollary~\ref{c:grad_div_jac}(ii) arguing as in Remark~\ref{r:aftermain}(ii).
\end{proof}

\begin{remark}
\label{r:generalop}
At this point it should be clear that Theorem~\ref{t:main1} and Corollaries~\ref{c:grad_div_jac} and \ref{c:grad_Lp} can be extended to more general differential operators. More precisely, let $T$ be an operator from $C^1_c(\R^d;\R^m)$ in $L^p(\mu;\R^n)$ of the form 
\begin{equation}
\label{e:generalop}
(Tu)_k := \sum_{j=1}^m T_{v_{jk}} u_j    
\end{equation}
where $v_{jk}$ are vector fields in $L^p(\mu)$ for every $1\le j\le m$ and $1\le k \le n$, 
and $T_{v_{jk}}$ is the directional derivative defined in \eqref{e:directional}. 
(Note that for $p=\infty$ this class includes all linear first order differential operator with constant coefficients.)

Then one can prove the following result:
if $v_{jk}(x) \in V(\mu,x)$ for $\mu$-a.e.~$x$ and every $j,k$, then $T$ can be extended to a continuous operator $\smash{\widetilde{T}} : \Lip(\R^d;\R^m) \to L^p(\mu;\R^n)$, namely the one obtained by replacing each $T_{v_{ij}}$ in formula \eqref{e:generalop} by the extension $\smash{ \widetilde{T}_{v_{ij}} }$ given in Theorem~\ref{t:main1}(i).
In particular $T$ is closable from $\Lip(\R^d)$ to $L^p_w(\mu)$.

Conversely, if $v_{jk}(x) \notin V(\mu,x)$ for a set of positive $\mu$-measure of points $x$ and for at least one couple of indices $j,k$, then $T$ is not closable from $\Lip(\R^d)$ to $L^p_w(\mu)$, nor from $\Lip(\R^d)$ to $L^p(\mu)$ with $p<\infty$, nor from $L^q(\mu)$ to $L^p(\mu)$. 
\end{remark}

\section{Closability of directional derivatives from $L^q$ to $L^p$}
\label{s:Lp}
\begin{theorem}
\label{t:lp_to_lp}
Let $1\le p,q\le\infty$ and let $p',q'$ denote the corresponding H\"older conjugates.
Let $v$ be a vector field in $L^p(\mu)$, and let $T_v: C^1_c(\R^d) \to L^p(\mu)$ be the directional derivative operator defined in \eqref{e:directional}. Assume that there exists a Borel function $\alpha$ on $\R^d$ such that
\begin{itemizeb}
\item
$\alpha\neq 0\; \mu$-a.e.;
\item
$\alpha\in L^{p'}(\mu)$ and $\alpha v\in L^{q'}(\mu)$;
\item
$N:=\alpha v\mu$ is a normal $1$-current.
\end{itemizeb}
Then $T_v$ is closable from $L^q(\mu)$ to $L^p(\mu)$.
\end{theorem}

\begin{remark}
\label{r:lp_to_lp}
Unlike Theorem~\ref{t:main1}(i), the closability result above is not accompanied by any differentiability result for the operator $\smash{\widetilde{T}_v}$ obtained by closing the graph of $T_v$ from $L^q$ to $L^p$. 

For example, let $\mu$ be the restriction of the Lebesgue measure to any bounded open set $\Omega$ in $\R^2$ and let $v(x):= (1,2x_1)$ and $\alpha(x):=1$ for every $x\in\R^2$. 
Then the assumptions in Theorem~\ref{t:lp_to_lp} are satisfied for every $1\le p,q \le \infty$, and $\smash{ \widetilde{T}_v } u = 0$ for every function $u:\R^2\to\R$ of the form $u(x_1,x_2):=g(x_2-x_1^2)$ with $g:\R\to\R$ bounded and Borel; moreover there exist plenty of $u$ of this form which are not differentiable at any point, in any direction.
\end{remark}


The rest of this section is devoted to the proof Theorem~\ref{t:lp_to_lp}.
Unspecified measures (as the expression $ds$ that appears in some integrals) are always the Lebesgue measure $\Leb^1$.
The proof relies on the next two lemmas.

\begin{lemma}
\label{l:currentdec}
Let $N=\widehat{\tau}\lambda$ be a normal $1$-current with $|\widehat{\tau}|=1$ $\lambda$-a.e.
Then $N$ can be decomposed as follows:
\begin{equation}
\label{e:currentdec}
N=\int_0^m[\gamma_s] \, ds
\end{equation}
where $m$ is a suitable positive number and
\begin{enumeraterem}
\item[\rm(a)]
for every $s\in [0,m]$, $\gamma_s:[0,L_s]\to\R^d$ is a Lipschitz path parametrized by arc-length, that is, $L_s$ is the length of $\gamma_s$ and $|\gamma_s'(t)|=1$ for a.e. $t \in [0,L_s]$;
\item[\rm(b)]
$|N|=\lambda=\int_0^m \lambda_s \, ds$ with $\lambda_s:= \big| [\gamma_s] \big|$;
\item[\rm(c)]
$\gamma'_s(t) = \widehat{\tau}(\gamma_s(t))$ for a.e.~$s\in [0,m]$ and a.e.~$t\in [0,L_s]$;
\item[\rm(d)]
$\lambda_s$ is the push-forward according to $\gamma_s$ of the Lebesgue measure on $[0,L_s]$.
\end{enumeraterem}
\end{lemma}

\begin{proof}
This decomposition is a variant of a well-know result by S.~Smirnov \cite{Smirnov}.
The result as stated can be found in \cite[Theorem~3.1]{PaoliniStepanov}, except for statement~(c), which can be proved as in \cite[Theorem~5.5(ii)]{AlbMar}, and statement~(d), which is a direct consequence of~(c) and the area formula. 
We refer to \cite{PaoliniStepanov} or \cite{AlbMar} for more details.%
\footnoteb{For example, the precise meaning of the integrals of in formula \eqref{e:currentdec} and statement~(b), and the correct measurability assumption on the map $s\mapsto[\gamma_s]$.}
\end{proof}

\begin{lemma}
\label{l:measuredec}
Let $\lambda$ be a finite positive measure on $\R^d$ that can be decomposed as 
$\lambda=\int_0^m \lambda_s\, ds$, and let $(v_n)$ be a sequence of functions such that $v_n \to v$ in $L^1(\lambda)$.
Then there exists a subsequence $(n_k)$ such that $v_{n_k} \to v$ in $L^1(\lambda_s)$ for a.e.~$s\in[0,m]$.
\end{lemma}

\begin{proof}
Let $g_n(s):=\int |v_n-v| \, d\lambda_s$.
Since $\int_0^m g_n(s) \, ds=\|v_n-v\|_{L^1(\lambda)}$, we have that $g_n \to 0$ in $L^1([0,m])$, and therefore there exists a subsequence $(n_k)$ such that $g_{n_k}(s) \to 0$ for a.e.~$s\in[0,m]$.
\end{proof}

\begin{proof}[Proof of Theorem~\ref{t:lp_to_lp}]
Since $T_v$ is linear, it suffices to prove that $T_v$ is ``closable at $0$'', namely that given a sequence $(u_n)\subset C^1_c(\R^d)$ and $w\in L^p(\mu)$ such that $u_n\to 0$ in $L^q(\mu)$ and $T_v u_n\to w$ in $L^p(\mu)$, then $w=0$.

Having set $E:=\{x\colon v(x)=0\}$, there holds $T_v u_n(x)=0$ for every $x\in E$ and every $n$, and therefore $w(x)=0$ for $\mu$-a.e.~$E$.
This means that it suffices to prove the statement above when $\mu$ is replaced by its restriction to $\R^d\setminus E$. 

In other words, we can freely assume that $v\ne 0$ $\mu$-a.e.

We set $\tau:=\alpha v$, and then $\tau\ne 0$ $\mu$-a.e.\ by the previous assumption.

We also set $\widehat{\tau}:=\tau/|\tau|$; then $|N|=|\tau|\mu$ and $N=\widehat{\tau}|N|$.

Since $u_{n}\to 0$ in $L^q(\mu)$ and $|\tau|\in \smash{L^{q'}}(\mu)$, then $u_n|\tau|\to 0$ in $L^1(\mu)$ or, equivalently, $u_n\to 0$ in $L^1(|\tau|\mu=|N|)$. 
Using the decomposition $|N|=\int_0^m \lambda_s\, ds$ in Lemma~\ref{l:currentdec}(b) and Lemma~\ref{l:measuredec} we have that, possibly passing to a subsequence in~$n$, $u_n\to 0$ in $L^1(\lambda_s)$ for a.e.~$s\in[0,m]$, and then, thanks to Lemma~\ref{l:currentdec}(d),
\[
u_n \circ \gamma_s\to 0
\quad\text{in $L^1([0,L_s])$ for a.e.~$s\in[0,m]$,}
\]
which in turn implies
\begin{equation}
\label{e:step1}
(u_n \circ \gamma_s)'\to 0
\quad\text{in $\D'(0,L_s)$ for a.e.~$s\in[0,m]$,}
\end{equation}
where $\D'(0,L_s)$ denotes the space of distributions on the interval $(0,L_s)$.

\smallskip
On the other hand, $T_v u_n =\nabla u_n \cdot v$ converges to $w$ in $L^p(\mu)$ by assumption, and since $\alpha\in\smash{L^{p'}}(\mu)$, then $\nabla u_n \cdot (\alpha v) \to \alpha w$ in $L^1(\mu)$.
Recalling that $\alpha v =\tau= \widehat{\tau} \, |\tau|$ we rewrite the last convergence as $(\nabla u_n \cdot \widehat{\tau}) |\tau| \to \alpha w$ in $L^1(\mu)$ or, equivalently, $\nabla u_n \cdot \widehat{\tau} \to \alpha w/|\tau|$ in $L^1(|\tau|\mu=|N|)$, and arguing as above we obtain that, possibly passing to a subsequence in~$n$, 
\begin{equation}
\label{e:step2}
(\nabla u_n\cdot \widehat{\tau}) \circ \gamma_s
\to 
(\alpha w/|\tau|) \circ \gamma_s
\quad\text{in $L^1([0,L_s])$ for a.e.~$s\in[0,m]$.}
\end{equation}
Thanks to Lemma~\ref{l:currentdec}(c) we obtain that, for every $n$ and a.e.~$s\in[0,m]$,
\[
(\nabla u_n\cdot \widehat{\tau}) \circ \gamma_s
= (\nabla u_n \circ \gamma_s) \cdot \gamma_s'
= ( u_n\circ \gamma_s)'
\quad\text{a.e.~in $[0,L_s]$,}
\]
and then \eqref{e:step2} becomes
\begin{equation}
\label{e:step3}
(u_n\circ \gamma_s)'
\to 
(\alpha w/|\tau|) \circ \gamma_s
\quad\text{in $L^1([0,L_s])$ for a.e.~$s\in[0,m]$}
\, .
\end{equation}
From \eqref{e:step1} and \eqref{e:step3} we infer that
\[
(\alpha w/|\tau|) \circ \gamma_s=0
\quad\text{a.e.~in $[0,L_s]$ for a.e.~$s\in[0,m]$.}
\]
By Lemma~\ref{l:currentdec}(d) this means that $\alpha w/|\tau|=0$ 
$\lambda_s$-a.e.\ and for a.e.~$s\in[0,m]$, which, by Lemma~\ref{l:currentdec}(b), implies
\[
\alpha w/|\tau|=0 
\quad\text{$|N|$-a.e., that is, $(|\tau|\mu)$-a.e.}
\]
Finally we use that $\alpha\ne 0$ and $\tau\ne 0$ $\mu$-a.e.\ to conclude that
$w=0$ $\mu$-a.e.
\end{proof}


	%
	%
\section{Closability of multilinear operators and metric currents}
\label{ss:metric}
In Theorem~\ref{t:multilin} we extend Theorem~\ref{t:main1} to the class of alternating $k$-linear differential operators defined in Paragraph~\ref{ss:klinear}.%
\footnoteb{For $k=n$ this class includes the Jacobian determinant; for that particular operator Theorem~\ref{t:multilin} reduces to Corollary~\ref{c:grad_div_jac}.}
We then study some closely related objects, namely metric $k$-currents in $\R^d$.
The key point is that metric $k$-currents can be naturally viewed as alternating $k$-linear differential operators (Remark~\ref{r-metriccurrents}(iii)).

\smallskip
Before stating Theorem~\ref{t:multilin} we recall the definition of $k$-tangent bundle of a measure and then define the class of $k$-linear operators we are interested in.

\begin{parag}[$\boldsymbol{k}$-tangent bundle]
\label{o:m-bundle}
(See \cite[Section~4]{AlbMar2} for more details.)
Given a measure $\mu$ on $\R^d$ and an integer $k$ with $1 \le k\le d$, the $k$-tangent bundle $V_k(\mu,\cdot)$ is a Borel map  $V_k(\mu,\cdot)$ on $\R^d$ whose values are vector subspaces of the space $\largewedge_k(\R^d)$ of $k$-vectors in $\R^d$ defined as follows: a $k$-vector $v$ belongs to $V_k(\mu,x)$ if and only if there exists a normal $k$-current $N$ in $\R^d$ with $\partial N=0$ such that \eqref{2.1} holds.%
\footnoteb{In other words, we have replaced the vector $v$ and the normal $1$-current $N$ in the definition of $V(\mu,x)$ by a $k$-vector and a normal $k$-current, respectively. In particular $V_1(\mu,x)=V(\mu,x)$.}

Given a $k$-vector $v$ in $\R^d$, we denote by $\Span(v)$ its supporting plane (or span), that is, the smallest subspace $W$ of $\R^d$ such that $v$ agrees with a $k$-vector on $W$ (cf.\ \cite{AlbMar}, Paragraph~5.8 and Proposition~5.9).

Note that if $v$ belongs to $V_k(\mu,x)$, then $\Span(v)$ is contained in $V(\mu,x)$ (\cite[Proposition~5.6]{AlbMar}). 
We do not know if the converse holds -- namely if every $k$-vector $v$ in $\R^d$ such that $\Span(v) \subset V(\mu,x)$ belongs to $V_k(\mu,x)$ -- except for the trivial cases $k=1$ and $k=d$ (see \cite{AlbMar2} for a more detailed discussion).
\end{parag}

\begin{parag}[Alternating $\boldsymbol{k}$-linear differential operators]
\label{ss:klinear}
Let $\tau$ be a $k$-vector field on $\R^d$ (that is, a map from $\R^d$ to $\largewedge_k(\R^d)$) which is bounded and Borel.
We denote by $J_\tau$ the alternating $k$-linear differential operator defined by
\begin{equation}
\label{e:klinear}
J_\tau u
:= \scalar{ \tau }{ du_1 \wedge\dots\wedge du_k }
\end{equation}
for every $u=(u_1,\dots,u_k)\in (C^1_c(\R^d))^k$.%
\footnoteb{Here $dv$ denotes the differential of the function $v$, intended as a $1$-form, and $\scalar{\,\cdot\,}{\cdot\,}$ is the duality pairing of $k$-vectors and $k$-covectors.}

In the next theorem we show that the closability and continuity properties of $J_\tau$ are connected to the following two assumptions on~$\tau$:
\begin{equation}
\label{e:multilin1}
\tau(x) \in V_k(\mu,x)
\quad\text{for $\mu$-a.e.~$x$,}
\end{equation}
which in turn implies (cf.~Paragraph~\ref{o:m-bundle})
\begin{equation}
\label{e:multilin2}
\Span(\tau(x)) \subset V(\mu,x)
\quad\text{for $\mu$-a.e.~$x$.}
\end{equation}
\end{parag}


\begin{theorem}
\label{t:multilin}
Take $\tau$ and $J_\tau$ be as above. 
\begin{enumeratethm}
\item
If \eqref{e:multilin2} holds, then every $u\in(\Lip(\R^d))^k$ is differentiable at $\mu$-a.e.~$x\in \R^d$ with respect to $W(x):=\Span(\tau(x))$, and the operator $\smash{ \widetilde{J}_\tau } : (\Lip(\R^d))^k\to L^\infty_{w}(\mu)$ given by
\begin{equation}
\label{e:multilin3}
\widetilde{J}_\tau u(x) 
:= \bigscalar{ \tau(x) }{ d_Wu_1(x) \wedge\dots\wedge d_Wu_k(x) }
\quad\text{for $\mu$-a.e.~$x$,}   
\end{equation}
is well defined and extends $J_\tau$. Moreover $\widetilde{J}_\tau (u_1,\dots,u_k)$ is separately continuous in each variable~$u_i$.

\item
Conversely, if $J_\tau (u_1,\dots,u_k)$ is separately continuous in each variable $u_i$ as a map from $C^1_c(\R^d)$ endowed with $\Lip$-convergence to $L^\infty_{w}(\mu)$, then \eqref{e:multilin2} holds.
\\
Accordingly, if \eqref{e:multilin2} does not hold, then  $J_\tau$ is not closable from $(\Lip(\R^d))^k$ to $L^\infty_w(\mu)$.

\item
If \eqref{e:multilin1} holds, then $\smash{ \widetilde{J}_\tau }$ is continuous;
it follows that $J_\tau$ is closable from $(\Lip(\R^d))^k$ to $L^\infty_w(\mu)$. 
\end{enumeratethm}
\end{theorem}

\begin{remark}
As pointed out in Paragraph~\ref{o:m-bundle}, we do not know if the sufficient condition for closability \eqref{e:multilin1} agrees with the necessary condition \eqref{e:multilin2}. Moreover we do not know if the necessary condition \eqref{e:multilin2} is also sufficient.
\end{remark}

\begin{proof}[Proof of Theorem~\ref{t:multilin}(i)]
Assumption \eqref{e:multilin2} implies the existence of the differential $d_Wu(x)$ for every $u\in\Lip(\R^d)$ and $\mu$-a.e.~$x$, and therefore $\smash{ \widetilde{J}_\tau }u$ is well defined. 

To prove the separate continuity of $\smash{ \widetilde{J}_\tau }$, we first rewrite it in terms of the directional derivative operators defined in  \eqref{e:extension}.
Given $\overline{u}_i\in \Lip(\R^d)$ for $i=1,\dots,k-1$, let $v$ be the ($1$-) vector field on $\R^d$ given by
\begin{equation}
\label{e:riduzione1}
v(x) := \tau(x) \trace 
    \big(
    d_W\overline{u}_1\wedge\cdots\wedge d_W\overline{u}_{k-1}
    \big)
\end{equation}
where $\trace$ is the interior product of $k$-vectors and $(k-1)$-covectors (in this specific case, vectors and covectors in the linear space $W(x)$), and define $\smash{ \widetilde{T}_v }$ according to~\eqref{e:extension}.
Then, for every $u\in\Lip(\R^d)$, 
\begin{align*}
\widetilde{J}_\tau ( \overline{u}_1, \dots, \overline{u}_{k-1}, u )
&   = \bigscalar{ \tau }{ d_W\overline{u}_1 \wedge \cdots \wedge d_W\overline{u}_{k-1} \wedge d_Wu } \\
&   = \bigscalar{ \tau \trace (d_W\overline{u}_1 \wedge \cdots \wedge d_W\overline{u}_{k-1}) }{ d_Wu } 
    = \scalar{ v }{ d_Wu }
    = T_v u 
    \, .
\end{align*}
Thus $\smash{\widetilde{J}}_\tau ( \overline{u}_1, \dots, \overline{u}_{k-1}, u )$ is continuous in $u$ by Theorem~\ref{t:main1}(i).
\end{proof}

\begin{proof}[Proof of Theorem~\ref{t:multilin}(ii)]
We exploit again the connection between the operator $J_\tau$ 
and the directional derivative operators defined in \eqref{e:directional}.
Given $1$-covectors $\alpha_i$ with $i=1,\dots,k-1$, let $v$ be the ($1$-) vector field on $\R^d$ given by
\begin{equation}
\label{e:riduzione}
v (x) := \tau(x) \trace (\alpha_1\wedge\cdots\wedge\alpha_{k-1})
\, ,
\end{equation}
and for every $i$ let $\overline{u}_i$ be the linear function on $\R^d$ such that $d\overline{u}_i=\alpha_i$.
Then for every $u\in C^1_c(\R^d)$ there holds
\begin{align*}
    J_\tau (\overline{u}_1, \dots, \overline{u}_{k-1}, u )
&   = \bigscalar{ \tau }{ \alpha_1\wedge\cdots\wedge\alpha_{k-1} \wedge du } \\
&   = \bigscalar{ \tau \trace (\alpha_1\wedge\cdots\wedge\alpha_{k-1}) }{ du }
    = \scalar{ v }{ du }
    = T_v u 
    \, .
\end{align*}
Since $J_\tau (\overline{u}_1, \dots, \overline{u}_{k-1}, u )$ is continuous in $u$, then $T_v$ is also continuous, and then Theorem~\ref{t:main1}(ii) implies $v(x)\in V(\mu,x)$ for $\mu$-a.e.~$x$.
Recalling \eqref{e:riduzione} and using the fact that every $(k-1)$-covector $\alpha$ is a linear combination of simple covectors of the form $\alpha_1\wedge\cdots\wedge\alpha_{k-1}$
we obtain that 
\[
\tau(x)\trace\alpha \in V(\mu,x)
\quad\text{for every $\alpha\in\largewedge^{k-1}(\R^d)$ and $\mu$-a.e.~$x$,}
\]
and now \eqref{e:multilin2} follows from the fact that $\Span(\tau(x))$ consists of all vectors of the form $\tau(x)\trace\alpha$ with $\alpha\in\largewedge^{k-1}(\R^d)$ (see \cite[Proposition~5.9]{AlbMar}).
\\
Finally, the non-closability of $J_\tau$ follows from the lack of continuity and the weak* pre-compactness of bounded subsets of $L^\infty(\mu)$.
\end{proof}

\begin{proof}[Proof of Theorem~\ref{t:multilin}(iii)]
Thanks to Theorems~1.1 and 1.2 in \cite{AlbMar2}, assumption \eqref{e:multilin1} implies that there exists a normal $k$-current $N=\tilde{\tau} \tilde{\mu}$ in $\R^d$ such that $\partial N=0$, and $\tilde{\tau}$ and $\tilde{\mu}$ extend $\tau$ and $\mu$ in the sense specified in the proof of Theorem~\ref{t:main1}(i), case $p=\infty$; proof that we follow almost verbatim to obtain the continuity of $\smash{ \widetilde{J}_\tau }$.
\end{proof}

We recall now the definition of metric currents; see \cite{metriccurrents}, \cite{lang} for more details.

\begin{parag}[Metric currents]
\label{ss:submetric}
Let $(X,d)$ be a complete metric space, and let $\Lip_b(X,\R)$ be the space of bounded Lipschitz functions on $X$. Given an $k\in\N$, a $k$-dimensional \emph{metric current} on $X$ is a functional $T:\Lip_b(X) \times (\Lip(X))^k \to\R$
that satisfies the following assumptions:
\begin{enumeratethm}
\item
\emph{linearity}: $T(f,\pi_1,\dots,\pi_k)$ is linear in each variable;
\item
\emph{continuity}: for every $f$, 
$T(f,\pi_1,\dots,\pi_k)$ is sequentially continuous in 
the variables $\pi_1,\dots,\pi_k$ \wrt pointwise convergence with uniformly bounded Lipschitz constants;
\item
\emph{locality}: $T(f,\pi_1,\dots,\pi_k)=0$ whenever there exists $i\in\{1,\dots,k\}$ such that $\pi_i$ is constant on a neighbourhood of $\supp(f)$;
\item
\emph{finite mass}: there exists a finite measure $\mu$ on $X$ such that, for every $f$ and every $\pi_i,\dots,\pi_k$,
\begin{equation}
\label{e:finitemass}
\big| T(f,\pi_1,\dots,\pi_k) \big|
\le \Lip(\pi_1)\cdots\Lip(\pi_k) \cdot \|f\|_{L^1(\mu)}
\, .
\end{equation}
\end{enumeratethm}
Finally, we define the support of $T$ as the smallest closed set $C\subset X$ such that $T(f,\pi_1,\dots,\pi_k)=0$ whenever $C\cap \supp(f) = \varnothing$.
\end{parag}

\begin{remarks}
\label{r-metriccurrents}
\begin{enumeraterem}
\item
Note that a metric current $T$ is \emph{alternating} in the variables $\pi_1,\dots,\pi_k$, that is, the value of $T(f,\pi_1,\dots,\pi_k)$ changes sign if we swap  $\pi_i$ and $\pi_j$ for any $i,j$ (this follows from the chain-rule in \cite[Theorem~3.5]{metriccurrents}).
\item 
If $T$ has compact support, the continuity assumption~(ii) is equivalent to say that $T(f,\pi_1,\dots,\pi_k)$ is sequentially continuous in the variables $\pi_1,\dots,\pi_k$ \wrt the usual $\Lip$-convergence.
\item 
Estimate \eqref{e:finitemass} implies that $T$ can be extended by continuity to all $f\in L^1(\mu)$, and therefore $T$ can be viewed as a ($k$-linear, alternating) operator that to every $(\pi_1,\dots,\pi_k) \in (\Lip(X))^k$ associate an element of $(L^1(\mu))'=L^\infty(\mu)$.
If $T$ has compact support, this operator is continuous from $(\Lip(X))^k$ to $L^\infty_w(\mu)$. 
\end{enumeraterem}
\end{remarks}

In the rest of this section $X$ is $\R^d$ endowed with the usual Euclidean distance. In this setting it is natural to compare metric currents and classical ones: the basic connection between these two notions is described in Paragraph~\ref{ss:metrictoclassical}; in Theorem~\ref{p:metric1} we give a new and more detailed description of such connection.

\begin{parag}[From metric to classical currents]
\label{ss:metrictoclassical}
To every metric $k$-current $T$ on $\R^d$ with compact support one can associate a classical $k$-current $\smash{ \widetilde{T} }$ defined as follows (cf.~\cite[Theorem~11.1]{metriccurrents}): 
\begin{equation}
\label{e:deftilda}
\bigscalar{ \widetilde{T} }{ \omega }
:= \sum_{\ii\in I(d,k)}
    T(\omega_\ii, x_{i_1},\dots, x_{i_k})
\end{equation}
for every $k$-form $\omega$ of class $C^1_c(\R^d)$, written in coordinates as
\[
\omega =\sum_{\ii\in I(d,k)}
\omega_\ii \, dx_{i_1}\wedge\dots\wedge dx_{i_k}
\, ,
\]
where $I(d,k)$ is the set of multi-indices $\ii=(i_1,\dots,i_k)$ with $1\leq i_1 < \dots < i_k\leq d$.
One easily checks that for every $f,\pi_1,\dots,\pi_K \in C^1_c(\R^d)$ there holds
\begin{equation}
\label{e:uguaglianza}
T(f,\pi_1,\dots,\pi_k) =
\bigscalar{ \widetilde{T} }{ f\, d\pi_1\wedge\cdots\wedge d\pi_k }
\, .
\end{equation}
\end{parag}

\begin{theorem}
\label{p:metric1}
The following statements hold:
\begin{enumeratethm}
\item
Let $T$ be a metric $k$-current on $\R^d$ with compact support.
Then there exists a finite measure $\mu$ with compact support and a bounded $k$-vector field $\tau$ on $\R^d$ such that \eqref{e:multilin2} holds, that is, $W(x):=\Span(\tau(x))$ is contained in $ V(\mu, x)$ for $\mu$-a.e.~$x\in\R^d$, and
\begin{equation}
\label{e:tmetric}
T(f,\pi_1,\dots,\pi_k)
=\int_{\R^d}
    f\, \scalar{\tau }{ d_W\pi_1\wedge\dots\wedge d_W\pi_k } \, d\mu  
\end{equation}
for all $f\in\Lip_b(\R^d)$ and all $\pi_1,\dots,\pi_k\in\Lip(\R^d)$.
\item 
On the other hand, given a finite measure $\mu$ with compact support and a bounded $k$-vector field $\tau$ such that \eqref{e:multilin1} holds, that is, $\tau(x) \in V_k(\mu, x)$ for $\mu$-a.e.~$x\in\R^d$, then formula \eqref{e:tmetric} defines a metric $k$-current $T$ with compact support.
\end{enumeratethm}
\end{theorem}

\begin{remarks}
\label{r:metric}
\begin{enumeraterem}
\item 
The assumption that $T$ and $\mu$ have compact support in statements~(i) and (ii) above is only needed to express the continuity assumption in the definition of metric currents (see Paragraph~\ref{ss:submetric}) in terms to the usual $\Lip$-convergence (cf.~Remark~\ref{r-metriccurrents}(ii)), and it can be removed with some care. 

\item 
A slight modification of the proof gives the following generalization of statement~(i): a functional $T : \Lip_b(\R^n) \times (\Lip(\R^d))^k \to \R$ admits an integral representation as in \eqref{e:tmetric} if and only if $T$ satisfies all assumptions in the definition of metric currents except continuity, which is replaced by the following weaker assumption: $T(f,\pi_1,\dots,\pi_k)$ is separately continuous in each variable $\pi_i$ \wrt $\Lip$-convergence.

\item
The assumption $\tau(x) \in V_k(\mu, x)$ for $\mu$-a.e.~$x\in\R^d$ in statement~(ii) is equivalent to say that $\tau\mu$ is a flat chain with finite mass (see Proposition~\ref{p:flatchains}(ii) below).
Thus formula \eqref{e:tmetric} defines a map from the space of flat chains with finite mass into metric $k$-currents,%
\footnoteb{This map has been already defined in \cite[Theorem~5.5]{lang}, but in a less explicit form.}
which is clearly a right inverse of the map $T\mapsto\smash{\widetilde{T}}$ defined in \eqref{e:deftilda}.
The Flat Chain Conjecture 
(see \cite[Section~11]{metriccurrents}, and Theorem~\ref{r:flat_conjecture}) states that this map is a bijection, that is, $\smash{\widetilde{T}}$ is a flat chain with finite mass for every metric current $T$.
\end{enumeraterem}
\end{remarks}

\begin{proof}[Proof of Theorem~\ref{p:metric1}(i)]
Let $\smash{ \widetilde{T} }$ be the (classical) current defined in \eqref{e:deftilda}.
Estimate \eqref{e:finitemass} yields
\[
\big| \bigscalar{ \widetilde{T} }{ \omega } \big| 
\le C \| \omega \|_{L^1(\mu)}
\]
where $C := \#(I(k,d)) = \binom{d}{k}$ and $\mu$ is the measure in Paragraph~\ref{ss:submetric}(iv).
Thus $\smash{ \widetilde{T} }$ has finite mass and can be written as $\smash{ \widetilde{T} }=\tau\mu$ where $\tau$ is a bounded $k$-vector field.
Recalling \eqref{e:uguaglianza} we obtain that for every $f,\pi_1,\dots,\pi_k \in C^1_c(\R^d)$ there holds
\begin{equation}
\label{e:uguaglianza2}
\begin{aligned}
    T(f,\pi_1,\dots,\pi_k)
&   = \bigscalar{ \widetilde{T} }{ f \, d\pi_1\wedge\cdots\wedge d\pi_k } \\
&   = \int_{\R^d} \scalar{\tau}{ d\pi_1\wedge\cdots\wedge d\pi_k } \, f \, d\mu 
    = \int_{\R^d} J_\tau(\pi_1,\dots\pi_k) \, f \, d\mu
    \, ,
\end{aligned}
\end{equation}
where $J_\tau: (C^1_c(\R^n))^k \to L^\infty(\mu)$ is defined in \eqref{e:klinear}.

Since $T$ is continuous in the variables $\pi_1,\dots,\pi_k$ \wrt $\Lip$-convergence (cf.~Paragraph~\ref{ss:submetric}(ii) and Remark~\ref{r-metriccurrents}(ii)), identity~\eqref{e:uguaglianza2} implies that $J_\tau$ is a continuous operator from $(C^1_c(\R^n))^k$, endowed with $\Lip$-convergence, to $L^\infty_w(\mu)$.

Then Theorem~\ref{t:multilin}(ii) implies that $\Span(\tau(x))\subset V(\mu, x)$ for $\mu$-a.e.~$x\in\R^d$, and therefore the integral at the right-hand side of \eqref{e:tmetric} is well defined.

To conclude we notice that identity \eqref{e:tmetric} holds for $f,\pi_1,\dots,\pi_k\in C^1_c(\R^d)$ by \eqref{e:uguaglianza2}, and can be extended 
to $f\in\Lip_b(\R^d)$ and $\pi_1,\dots,\pi_k\in\Lip(\R^d)$ by continuity.

Indeed, the continuity of the left-hand side of \eqref{e:tmetric} follows from the definition of metric currents, while the separate continuity of the right-hand side \wrt each of the variables $f,\pi_1,\dots,\pi_k$ follows from the separate continuity of the operator $\smash{ \widetilde{J} }_\tau$ proved in Theorem~\ref{t:multilin}(i).
\end{proof}

\begin{proof}[Proof of Theorem~\ref{p:metric1}(i)]
Formula~\eqref{e:tmetric} can be re-written as
\begin{equation}
\label{e:T=J}
T(f,\pi_1,\dots,\pi_k)
=\int_{\R^d}
    \widetilde{J}_\tau(\pi_1,\dots,\pi_k) \, f\, d\mu  
\end{equation}
where $\smash{\widetilde{J}}_\tau$ is defined as in \eqref{e:multilin3}.
Then $T$ satisfies the continuity assumption in the definition of metric currents because $\smash{\widetilde{J}}_\tau$ is continuous from $(\Lip(\R^n))^k$ to $L^\infty_w(\mu)$ (Theorem~\ref{t:multilin}(iii)). 
It is immediate to check that $T$ satisfies all other assumptions in the definition of metric currents. 
\end{proof}

%
%


We conclude this section by pointing out that the Flat Chain Conjecture is equivalent to the converse of Theorem~\ref{t:multilin}(iii) and Theorem~\ref{p:metric1}(ii). 
Before giving a precise statement we recall the definition of flat chains and some useful characterizations.

\begin{parag*}[Flat chains]
The space of $k$-dimensional flat chains in $\R^d$ is defined as the closure of $k$-normal currents with respect to the flat norm.
The next statement gives some characterizations of flat chains \emph{with finite mass}.
\end{parag*}

\begin{proposition}
\label{p:flatchains}
Let $T=\tau\mu$ be a $k$-current with finite mass in $\R^d$, where $\mu$ is a finite measure and $\tau$ a $k$-vector field such that $\tau\ne 0$~$\mu$-a.e.
\begin{enumeratethm}
\item
If $k=d$, $T$ is a flat chain if and only if $\mu\ll\Leb^d$.
\item 
If $k\le d$, $T$ is a flat chain if and only if \eqref{e:multilin1} holds.
\item 
if $k <d$, $T$ is a flat chain if and only if it can be written as the restriction of a normal $k$-current to a Borel set.
\end{enumeratethm}
\end{proposition}

\begin{proof}
The proof of statement (i) is immediate.
Statement~(ii) follows from  \cite[Theorem~1.2]{AlbMar2}.
Statement~(iii) follows from  \cite[Theorem~1.1]{AlbMar2}.
\end{proof}

\begin{theorem}
\label{r:flat_conjecture}
Let $k=1,\dots,d$. The following statements are equivalent:
\begin{enumeratethm}
\item
\emph{(Flat Chain Conjecture)} 
Let $T$ be a metric $k$-current. Then the current $\smash{\widetilde{T}}$ defined in \eqref{e:deftilda} is a flat chain.
\item 
\emph{(Converse of Theorem~\ref{p:metric1}(ii))}
Let $T$ be a metric $k$-current, and let $\mu$ and $\tau$ be as in 
Theorem~\ref{p:metric1}(i). Then \eqref{e:multilin1} holds.
\item 
\emph{(Converse of Theorem~\ref{t:multilin}(iii))}
Let $\mu$ be a finite measure, $\tau$ a bounded $k$-vector field, and $J_\tau$ the operator defined in \eqref{e:klinear}. If $J_\tau$ is closable from $(\Lip(\R^d))^k$ to $L^\infty_w(\mu)$ then \eqref{e:multilin1} holds.
\end{enumeratethm}
\end{theorem}

\begin{proof}
The equivalence of statements~(i) and (ii) is an immediate consequence of
Theorem~\ref{p:metric1} and Proposition~\ref{p:flatchains}(ii).
The equivalence of statements~(ii) and (iii) is an immediate consequence of
the following facts: $J_\tau$ is closable if and only if $\smash{\widetilde{J}}_\tau$ is well-defined and continuous (Theorem~\ref{t:multilin}), and the continuity of $\smash{\widetilde{J}}_\tau$ is equivalent to that of $T$ because of identity~\eqref{e:T=J}.
\end{proof}

The following corollary is well known (see \cite[Theorem~1.6]{Schioppa} and \cite[Theorem~1.15]{DPR}, and the more recent proofs \cite{DeMaMa}, \cite{MaMe});
we simply remark that it follows from our previous results as well. 

\begin{corollary}
\label{c:metric}
The Flat Chain Conjecture is true for $k=1$ and $k=d$.
\end{corollary}

\begin{proof}
For $k=1$ and $k=d$ properties \eqref{e:multilin1} and \eqref{e:multilin2} are equivalent (the case $k=1$ is trivial, while the case $k=d$ follows from Theorem~\ref{c:bundlesingular}).
This means that Theorem~\ref{p:metric1}(i) implies statement~(ii) in Theorem~\ref{r:flat_conjecture}, which implies the Flat Chain Conjecture by the very same theorem.
\end{proof}

    %
    %
    %
    %
\bibliography{ref}

\providecommand{\bysame}{\leavevmode\hbox to3em{\hrulefill}\thinspace}
\providecommand{\MR}{\relax\ifhmode\unskip\space\fi MR }
\providecommand{\MRhref}[2]{%
  \href{http://www.ams.org/mathscinet-getitem?mr=#1}{#2}
}
\providecommand{\href}[2]{#2}
\begin{thebibliography}{10}

\bibitem{ACP2}
Giovanni Alberti, Marianna Cs\"{o}rnyei, and David Preiss,
  \emph{Differentiability of {L}ipschitz functions, structure of null sets, and
  other problems}, Proceedings of the {I}nternational {C}ongress of
  {M}athematicians. {V}olume {III}, Hindustan Book Agency, New Delhi, 2010,
  pp.~1379--1394.

\bibitem{AlbMar}
Giovanni Alberti and Andrea Marchese, \emph{On the differentiability of
  {L}ipschitz functions with respect to measures in the {E}uclidean space},
  Geom. Funct. Anal. \textbf{26} (2016), no.~1, 1--66,
  doi:~\href{https://doi.org/10.1007/s00039-016-0354-y}{10.1007/s00039-016-0354-y}.

\bibitem{AlbMar2}
\bysame, \emph{On the structure of flat chains with finite mass}, {a}r{X}iv
  preprint \href{https://arxiv.org/abs/2311.06099}{2311.06099} (2023).

\bibitem{Albev}
Sergio Albeverio and Michael R\"ockner, \emph{Classical {D}irichlet forms on
  topological vector spaces -- {C}losability and a {C}ameron-{M}artin formula},
  J. Funct. Anal. \textbf{88} (1990), no.~2, 395--436,
  doi:~\href{https://doi.org/10.1016/0022-1236(90)90113-Y}{10.1016/0022-1236(90)90113-Y}.

\bibitem{AmbDal}
Luigi Ambrosio and Gianni Dal~Maso, \emph{A general chain rule for
  distributional derivatives}, Proc. Amer. Math. Soc. \textbf{108} (1990),
  no.~3, 691--702,
  doi:~\href{https://doi.org/10.2307/2047789}{10.2307/2047789}.

\bibitem{metriccurrents}
Luigi Ambrosio and Bernd Kirchheim, \emph{Currents in metric spaces}, Acta
  Math. \textbf{185} (2000), no.~1, 1--80,
  doi:~\href{https://doi.org/10.1007/BF02392711}{10.1007/BF02392711}.

\bibitem{Boga}
Vladimir~I. Bogachev, \emph{Differentiable measures and the {M}alliavin
  calculus}, Mathematical Surveys and Monographs, vol. 164, American
  Mathematical Society (AMS), Providence, RI, 2010.

\bibitem{BouButSep}
Guy Bouchitt\'e, Giuseppe Buttazzo, and Pierre Seppecher, \emph{Energies with
  respect to a measure and applications to low dimensional structures}, Calc.
  Var. Partial Differential Equations \textbf{5} (1997), 37--54,
  doi:~\href{https://doi.org/10.1007/s005260050058}{10.1007/s005260050058}.

\bibitem{BouChaJim}
Guy Bouchitt\'e, Thierry Champion, and Chlo\'e Jimenez, \emph{Completion of the
  space of measures in the {K}antorovich norm}, Riv. Mat. Univ. Parma (7)
  \textbf{4*} (2005), 127--139.

\bibitem{BreGig}
Camillo Brena and Nicola Gigli, \emph{About the general chain rule for
  functions of bounded variations}, Nonlinear Anal. \textbf{242} (2024),
  113518,
  doi:~\href{https://doi.org/10.1016/j.na.2024.113518}{10.1016/j.na.2024.113518}.

\bibitem{DeMaMa}
Luigi De~Masi and Andrea Marchese, \emph{A refined {L}usin type theorem for
  gradients}, {a}r{X}iv preprint
  \href{https://arxiv.org/abs/2411.15012}{2411.15012} (2024).

\bibitem{DPR}
Guido De~Philippis and Filip Rindler, \emph{On the structure of {$\mathcal
  A$}-free measures and applications}, Ann. of Math. (2) \textbf{184} (2016),
  no.~3, 1017--1039,
  doi:~\href{https://doi.org/10.4007/annals.2016.184.3.10}{10.4007/annals.2016.184.3.10}.

\bibitem{DLP}
Simone Di~Marino, Danka Lu{\v{c}}i{\'c}, and Enrico Pasqualetto, \emph{A short
  proof of the infinitesimal {H}ilbertianity of the weighted {E}uclidean
  space}, C. R., Math., Acad. Sci. Paris \textbf{358} (2020), no.~7, 817--825,
  doi:~\href{https://doi.org/10.5802/crmath.88}{10.5802/crmath.88}.

\bibitem{Federer1996GeometricTheory}
Herbert Federer, \emph{{G}eometric {M}easure {T}heory}, Grundlehren der
  Mathematischen Wissenschaften, vol. 153, Springer-Verlag,
  Berlin-Heidelberg-New York, 1969, reprinted in the series Classics in
  Mathematics, Springer, Berlin 1996,
  doi:~\href{https://doi.org/10.1007/978-3-642-62010-2}{10.1007/978-3-642-62010-2}.

\bibitem{Fuk}
Masatoshi Fukushima, Yoichi Oshima, and Masayoshi Takeda, \emph{Dirichlet forms
  and symmetric {M}arkov processes}, second ed., De~Gruyter Studies in
  Mathematics, vol.~19, De~Gruyter, Berlin, 2010,
  doi:~\href{https://doi.org/10.1515/9783110889741}{10.1515/9783110889741}.

\bibitem{KrantzParks}
Steven~G. Krantz and Harold~R. Parks, \emph{Geometric integration theory},
  Cornerstones, Birkh\"auser, Boston-Basel-Berlin, 2008,
  doi:~\href{https://doi.org/10.1007/978-0-8176-4679-0}{10.1007/978-0-8176-4679-0}.

\bibitem{lang}
Urs Lang, \emph{Local currents in metric spaces}, J. Geom. Anal. \textbf{21}
  (2011), no.~3, 683--742,
  doi:~\href{https://doi.org/10.1007/s12220-010-9164-x}{10.1007/s12220-010-9164-x}.

\bibitem{Lindenstrauss2012}
Joram Lindenstrauss, David Preiss, and Jaroslav Ti\v{s}er, \emph{{F}r\'echet
  differentiability of {L}ipschitz functions and porous sets in {B}anach
  spaces}, Annals of Mathematics Studies, vol. 179, Princeton University Press,
  Princeton, NJ, 2012,
  doi:~\href{https://doi.org/10.1515/9781400842698}{10.1515/9781400842698}.

\bibitem{MaMe}
Andrea Marchese and Andrea Merlo, \emph{A simple proof of the 1-dimensional
  flat chain conjecture}, {a}r{X}iv preprint
  \href{https://arxiv.org/abs/2411.15019}{2411.15019} (2024).

\bibitem{PaoliniStepanov}
Emanuele Paolini and Eugene Stepanov, \emph{Structure of metric cycles and
  normal one-dimensional currents}, J. Funct. Anal. \textbf{264} (2013), no.~6,
  1269--1295,
  doi:~\href{https://doi.org/10.1016/j.jfa.2012.12.007}{10.1016/j.jfa.2012.12.007}.

\bibitem{RenRog}
Michael Renardy and Robert~C. Rogers, \emph{An introduction to partial
  differential equations}, second ed., Texts in Applied Mathematics, vol.~13,
  Springer, New York, 2004,
  doi:~\href{https://doi.org/10.1007/b97427}{10.1007/b97427}.

\bibitem{Schioppa}
Andrea Schioppa, \emph{Metric currents and {A}lberti representations}, J.
  Funct. Anal. \textbf{271} (2016), no.~11, 3007--3081,
  doi:~\href{https://doi.org/10.1016/j.jfa.2016.08.022}{10.1016/j.jfa.2016.08.022}.

\bibitem{Smirnov}
Stanislav~K. Smirnov, \emph{Decomposition of solenoidal vector charges into
  elementary solenoids and the structure of normal one-dimensional currents},
  St. Petersbg. Math. J. \textbf{5} (1994), no.~4, 841--867, translation from
  Algebra i Analiz 5 (1993), 206--238.

\end{thebibliography}

\bibliographystyle{amsplain}

	%
	%
	%
	%
\vskip .5 cm

{\parindent = 0 pt\footnotesize
G.A.
\\
Dipartimento di Matematica, 
Universit\`a di Pisa,
largo Pontecorvo~5, 
56127 Pisa, 
Italy 	
\\
e-mail: \texttt{giovanni.alberti@unipi.it}

\bigskip
D.B.
\\
Mathematics Institute, 
Zeeman Building, 
University of Warwick, 
Coventry CV4 7AL, 
UK
\\
e-mail: \texttt{david.bate@warwick.ac.uk}

\bigskip
A.M.
\\
Dipartimento di Matematica, 
Universit\`a di Trento,
via Sommarive 14, 
38123 Povo, 
Italy 
\\
e-mail: \texttt{andrea.marchese@unitn.it}

}

\end{document}